\documentclass[12pt,a4paper]{article}
\pdfoutput=1
\usepackage{authblk}
\usepackage{graphicx}
\usepackage{caption}
\usepackage{float}
\usepackage{epstopdf}
\usepackage{amsmath,array,amsthm}
\usepackage{hyperref}
\usepackage{cleveref}
\usepackage{amssymb}
\usepackage{bm}
\usepackage{color}
\usepackage{geometry}
\usepackage{longtable}
\usepackage{multirow}
\usepackage{makecell}
\usepackage{enumerate}
\usepackage{cite}
\newtheorem{thm}{Theorem}
\newtheorem{ass}{Assumption}
\newtheorem{lem}{Lemma}
\newtheorem{defi}{Definition}
\newtheorem{rem}{Remark}
\usepackage{algorithm,algorithmic}
\usepackage{indentfirst}
\newcommand{\X}{\mathcal{X}}
\newcommand{\R}{\mathbb{R}}
\newcommand{\K}{\mathbb{K}}
\newcommand{\N}{\mathcal{N}}
\newcommand{\A}{\mathcal{A}}
\newcommand{\B}{\mathcal{B}}
\newcommand{\rc}{\mathcal{R}}
\newcommand{\mL}{\mathrm{L}}
\newcommand{\cL}{\mathcal{L}}
\title{Push-Pull Based Distributed Primal-Dual Algorithm  for Coupled Constrained Convex Optimization in Multi-Agent Networks}
\author[1]{Kai Gong }
\author[1]{Liwei Zhang}
\affil[1]{School of Mathematical Sciences, Dalian University of Technology, Dalian,
	116024, Liaoning, China. \authorcr Email: gk1995\_\_\_@mail.dlut.edu.cn
	\authorcr Email: lwzhang@dlut.edu.cn}
\date{}
\begin{document}
	\maketitle
	\begin{abstract}
		This paper focuses on a distributed coupled constrained convex optimization problem over directed unblanced and  time-varying multi-agent networks, where the global objective function is the sum of all agents' private local objective functions, and decisions of all agents are subject to coupled equality and inequality	constraints and a compact convex subset. In the multi-agent networks, each agent exchanges information with other neighboring agents. Finally, all agents reach a consensus on decisions, meanwhile achieve the goal of minimizing the global objective function under the given constraint conditions. For the purpose of protecting the information privacy of each agent, we first establish the saddle point problem of the constrained convex optimization problem considered in this article, then based on the push-pull method, develop a distributed primal-dual algorithm to solve the dual problem. Under the Slater's condition, we will show that the sequence of points generated by the proposed algorithm converges to a saddle  point of the Lagrange function. Moreover, we analyze the iteration complexity of the algorithm. \\
		
		\textbf{Keyworkds}: Distributed optimization algorithm, coupled constraints, time-varying networks, convergence, iteration complexity.
		
	\end{abstract}

	\section{Introduction}\label{sec:1}
	We consider the following distributed constrained optimization problem in directed and time-varying networks consisting of $m$ agents:
	\begin{equation}\label{P}
		\begin{aligned}
			\min_{x\in\X}&f(x)=\sum_{i=1}^mf_i(x)\\
			{\rm s.t.} &\quad \sum_{i=1}^mh_i(x)\in\K,\\
		\end{aligned}
	\end{equation}
	where $\K\triangleq\R^p_-\times\{0_q\}$, implying that the constraints include inequality and equality. Functions $f_i:\R^n\to \R$ and $h_i:\R^n\to \R^{p+q}$ are the local objective function and local constraint function of agent $i$, respectively, $i=1,2,\dots,m$, and $\X$ is a subset of $\R^n$. The local information $f_i$ and $h_i$ are only known by agent $i$, but other agents don't have direct access to them.
	
	\par Distributed convex optimization has acctracted much attention in recent years, due to broad range of applications in such fields as machine learning \cite{Mahmoudi20}, wireless sensor networks \cite{Liao21b}, social networks \cite{Sasso20}, smart grids \cite{Sharma16}, etc. For the distributed convex optimization problem  considered in this paper, the goal is to minimize a global objective function, given as a sum of local objective functions, subject to some constraints, which include coupled equality and ineqaulity constraints and a common constraint set. In this system, each agent has its own decision, and communicates with neiborning agents to adjust its decision individually. All agents expect to reach a agreement on decision which is an optimal solution to the considered optimization problem. The refered distributed convex optimization problem is motivated by its pratical applications in such fields as distributed estimation \cite{Nowak03}, distributed source localization \cite{Feng13}, network utility maximization \cite{Liao21a}, optimal flow control in power systems \cite{Shaoqun16}, etc.
	\subsection{Realted Work}
	For distributed problems, the unconstrained distributed convex optimization models were first widely studied. A seminal work in this direction  is the distributed  first-order (sub)gradient descent algorithm (DGD) proposed by \cite{Nedic09a}. Paper \cite{Nedic09a} established the asympotical convergence of DGD algorithm with diminishing stepsizes. For the case where the local objective functions are strongly convex and L-smooth, \cite{Shi15a} developed an exact first-order algorithm (abbreviated as EXTRA). In addition, \cite{Nedic17a} developed a distributed inexact gradient method with the gradient tracking technique (abbreviated as DIGing), which was also developed independently in \cite{Qu18, Xu15, Xu18}. The two algorithms, EXTRA and DIGing both achieved linear convergence rates for strongly convex and L-smooth local objective functions. While the weight matrices of EXTRA are fixed, and that of DIGing are doubly stochastic. During information exchange among agents, it is difficult to keep the adjacency matrices being doubly stochastic all the time. The push-pull method proposed by \cite{Pu18} could overcome this challenge, owing to that the method utilizes two different graphs for the information exchange among agents. Furthermore, the TV-AB algorithm proposed by \cite{Saad20} extended the push-pull method to the time-vary and directed unblanced networks. The graph at each iteration in TV-AB algorithm can be generated randomly in an arbitrary way as long as the union of eyery B consecutive graphs is strongly connected, which overcomes the difficulty of keeping the adjacency matrices being doubly stochastic at all times during information exchanging among all agents. Besides, there are many other push-pull based algorithms to solve different types of distributed problems, e.g., see \cite{Du18, Xin18, Xin19}.
	
	In real world, the decisions of agents located in multi-agent networks are usually subject to certain constraints. Naturally, the distributed constrained optimization problems  have  been  gradually devoted more efforts and interests to study. The earlier studies in this direction focused on the special case where the local constraint sets of agents are identical (common). For instance, \cite{Nedic10} developed a distributed projected (sub)gradient algorithm, and established the asymptotical convergence of the proposed algorithm. When the local constraint sets are different, \cite{Lei16} developed a distributed primal-dual algorithm, which can be viewed as an extension of the EXTRA to the constrained case. Besides, \cite{Zhu12} extended the distributed projected (sub)gradient method to a more general case with global (common) equlity and  inequality constraints and a common simple constraint set. 
	
	\par The distributed optimization problems with coupled inequality and equality constraints are also extremely important, which attracted much attention and interests from researchers in the related fields. Currently, only a few works have studied such problems. \cite{Carli20} proposed a distributed ADMM algorithm to handle the case with coupled affine equality constraints. In  the case with coupled inequality constraints, \cite{Falsone17} applied the combination of dual decomposition method and the distributed projection (sub)gradient method. Also, for the same problem, \cite{Chang14} adopted a consensus-based primal-dual perturbatin method in directed balanced and time-varying networks. When the considered network is static and undirected, \cite{Liang20} developed a distributed algorithm for the case with coupled inequality and affine equality constraints. Other related works can be found in \cite{Yang19, Zhang20, Nedic18}, et al.
	\subsection{Main Contribution}
	The main contribution of the article is threefold. First, the coupled constraints of the optimization problem considered in this article are more general than \cite{Falsone17, Chang14, Yang19, Zhang20, Nedic18}, which include coupled inequality (not necessarily affine) and affine equality constraints. Second, we develop a push-pull based distributed primal-dual algorithm to solve the challenging problem, and analyze the iteration complexity of the presented algorithm. Third, the considered multi-agent networks are more general than \cite{Falsone17, Chang14, Liang20, Zhang20, Nedic18}, which are directed weight-unblanced and time-varying. The requirement of keeping the adjacency matrices being doubly stochastic at all times during  information exchanging among agents is difficult, but our algorithm overcomes this difficulty, owing to the push-pull method. 
	
	\subsection{Organization of This Paper}
	The rest of this paper is organized as follows. \Cref{sec:2} provides necessary preliminaries, constructs the dual problem and develops a push-pull based distributed primal-dual algorithm. \Cref{sec:3} presents the main results about properties and convergence of the proposed algorithm. \Cref{sec:4} simulates the given algorithm for an example. Finally, \Cref{sec:5} makes some conclusions. 
	
	\section{Algorithm}\label{sec:2}    
	\subsection{Notations and Notion}
	We only consider the case where the multi-agent networks operate synchronously. The topology of the network at time $k\ge0$ will be represented by a directed weighted graph $G(k)=(V,E(k))$, where $V=\{1,2,\dots,m\}$ is the collection of all agents, and $E(k)$ is the set of directed edges. $A(k)$ and $B(k)$ of $n\times n$ matrices are compliant with $G(k)$,  i.e., $(i,j)\in E(k)\Leftrightarrow a_{ij}(k)>0$ and $b_{ij}(k)>0$. All agents which transmit information to node $i$ directly at time $k$ are said to be its in-neighbors and belong to the set $\N_i^{\rm{in}}(k)=\{j\in V: (j,i)\in E(k)\}$. All agents which receive information from agent $i$ at time $k$ belong to the set of its out-neighbors, denoted by $\N_i^{\rm{out}}(k)=\{j\in V: (i,j)\in E(k)\}$. In convergence analysis of the presented algorithm, we will use the nations of sochastic vector, absolute probability sequence and saddle point of a convex-concave function.
	\begin{defi}[Stochastic Vector]
		A vector $\pi\in\R^m$ is said to be a stochastic vector when its components $\pi_i, i=1,2,\dots,m$ satisfy that
		\begin{equation*}
			\sum_{i=1}^m\pi_i=1,\quad\pi_i\ge0,\quad \forall i=1,2,\dots,m. 
		\end{equation*}
	\end{defi}
	\begin{defi}[Absolute Probability Sequence]
		For row-stochastic matrices $\{\rc_k\}$, an absolute probability sequence is a sequence $\{\pi_k\}$ of stochastic vectors such that
		\begin{equation*}
			\pi_k^\top=\pi_{k+1}^\top\rc_k,\quad \forall k\ge0.
		\end{equation*}
	\end{defi}
	\begin{defi}[Saddle Point] 
		Consider a function $\cL:X\times S\to \R$, where $X$ and $S$ are non-empty subsets of $\R^n$ and $\R^{p+q}$ respectively, a pair of vectors $(x^*,\lambda^*)\in X\times S$ is called a saddle point of function $\cL$ over $ X\times S$, if \ $\forall(x,\lambda)\in X\times S$, it holds that
		\begin{equation*}
			\cL(x^*,\lambda)\le \cL(x^*,\lambda^*)\le \cL(x,\lambda^*).
		\end{equation*}	   
	\end{defi}
	In this paper, the gradient vector of a function $f$ at point $x$ is denoted by $\nabla f(x)$, $\Vert\cdot\Vert$ represents the Frobenius norm of a matrix of suitalbe dimensions, $\langle \cdot,\cdot\rangle$ denotes the inner product matched with the norm $\Vert\cdot\Vert$. The projection of a point $x$ onto a subset $\Omega\subseteq \R^n$ is  denoted by $P_\Omega(x)$. For a non-empty subset $\Omega$, we denote the relative interior of $\Omega$ by $\rm{ri}(\Omega)$. $1_m$ denotes the column vector in $\R^m$ whose entries are all ones.
	
	\par We make the following assumptions on the network communication graphs, which are standard  in the analysis of average consensus algorithms, e.g., see \cite{Nedic09a, Nedic10, Saad20}. 
	\begin{ass}[Periodical Strong Connectivity]\label{ass1}
		There exists a positive ineger $B$ such that $\forall k\ge0$, the directed graph $\left(V,\bigcup_{t=0}^{B-1}E(t+k)\right)$ is strongly connected.
	\end{ass}
	
	\begin{ass}[Weights]\label{ass2}
		\hspace{1em}
		\begin{enumerate}[\hspace{1em}(1)]
			\item Stochasticity: Matrices $\{A(k)\}$ and $\{B(k)\}$ are row-stochastic and column-stochastic, respectively.
			\item Aperiodicity: Graph $G(k)$ has self-loops, i.e., $a_{ii}(k)>0, b_{ii}(k)>0, \forall i\in V, k\ge 0$. 
			\item Uniform positicity: There is a scalar $\eta\in (0,1)$ such that $a_{ij}(k)\ge\eta$ and $b_{ij}(k)\ge\eta, \forall (i,j)\in E(k), k\ge0$.
		\end{enumerate}	 
	\end{ass}
	\subsection{Saddle point Problem}
	Smooth convex optimization problems are within the range of our consideration, then we make the following standard assumptions for the analysis of the formulated algorithm.
	\begin{ass}[Convexity and Smoothness]\label{ass3}
		\hspace{1em}
		\begin{enumerate}[\hspace{1em}(1)]
			\item The local function $f_i$ of agent $i$ is convex and continuous differentiable over the subset $\X$, $\forall i=1,2,\dots,m$.
			\item The coordinate component function $h_i^j$ of $h_i$ corresponding to the inequality constraints is convex and continuous differentiable, $j=1,\dots,p, i=1,\dots,m$. Function $h_i^j$ corresponding to the equality constraints is affine, i.e., $h_i^j(x)=(a_i^j)^\top x+b_{ij}, j=p+1,\dots,p+q$, where $a_i^j$ is the $j$-\rm{th} column of a matrix $A_i\in\R^{n\times p}$, $b_{ij}\in\R$, $\forall i=1,\dots,m$.
			\item The subset $\X$ of $\R^n$ is compact and convex.
		\end{enumerate}	
	\end{ass}
	
	\begin{ass}[Slater's Condition]\label{ass4}
		Suppose that there exists a point $\bar{x}\in \rm{ri}(\X)$, such that $\sum_{i=1}^mh_i(\bar{x})\in\rm{ri}(\K)$.
	\end{ass}
	Let $\K^\circ$ represent the polar of closed convex cone $\K$, then $\K^\circ=\R^p_+\times\R^q$. The dual problem of the constrained convex  optimization problem \eqref{P} can be formulated as follows:
	\begin{equation}\label{D}
		\max_{\lambda\in\K^\circ}\min_{x\in\X}\cL(x,\lambda)
	\end{equation}
	where $\cL$ is the Lagrange function defined by 
	\begin{align*}
		\cL(x,\lambda)&=\sum_{i=1}^m\cL_i(x,\lambda)=\sum_{i=1}^m\left\{f_i(x)+\lambda^\top h_i(x)\right\},\\ \lambda&=(\lambda^I,\lambda^E)\in\K^\circ,
	\end{align*} 
	here $\lambda^I$ and $\lambda^E$ are the multiplies corresponding to the inequality constraints and equality constraints respectively.
	\par Under Slater's condition, the strong duality threorem holds, i.e., there doesn't exist duality gap between primal problem \eqref{P} and dual problem \eqref{D}. In other words, a point $x^*\in\X$ is a solution to primal problem \eqref{P} if and only if there is a pair of points $(x^*,\lambda^*)\in\X\times\K^\circ$ such that $(x^*,\lambda^*)$ is a saddle point of the Lagrange function $\cL(x,\lambda)$. For any $\bar{\lambda}\in\K^\circ$, let
	\begin{align*}
		q(\lambda)&=\inf_{x\in\X}\cL(x,\lambda),\\
		Q(\bar{\lambda})&=\left\{\lambda=(\lambda^I,\lambda^E)\in\K^\circ:q(\lambda)\ge q(\bar{\lambda})\right\}.  
	\end{align*}
	The boundness of multiplier $\lambda^I$ corresponding to the inequality constraints is indispensable in the proof of convergence of our algorithm, which can be deduced from the following lemma.
	\begin{lem}[\cite{Nedic09b} Lemma1]\label{lem:1}
		For any $\bar{\lambda}\in\K^\circ,\alpha\in\R$, it holds that
		\begin{equation}
			\max_{\lambda\in Q(\bar{\lambda})}\Vert\lambda^I\Vert\le\frac{1}{\gamma(\bar{x})}\left(f(\bar{x})-q(\bar{\lambda})\right).
		\end{equation}
		where $\gamma(\bar{x})=\min_{1\le j\le p}\left\{-\sum_{i=1}^mh_i^j(\bar{x})\right\}$, $\bar{x}$ satisfies Slater's condition.
	\end{lem}
	\par Let $q^*$ denote the dual optimal value, and $Q^*$ denote the dual optimal solutions set, namely 
	$$Q^*\triangleq\left\{\lambda=(\lambda^I,\lambda^E)\in\K^\circ:q(\lambda)\ge q^*\right\}\subseteq Q(\bar{\lambda}).$$
	According to the above lemma, we have that the dual optimal multiplier $\lambda^I$ is bounded. Select an arbitrary vector $\bar{\lambda}\in\K^\circ$, we define 
	\begin{align*} 
		Q_I&\triangleq\left\{\lambda^I\in\R^p_+: \Vert\lambda^I\Vert\le\frac{1}{\gamma(\bar{x})}\left(f(\bar{x})-q(\bar{\lambda})\right)\right\},\\
		Q&\triangleq Q_I\times \R^q.
	\end{align*}
	Then we have $Q^*\subseteq Q(\bar{\lambda})\subseteq Q$. Next we will construct a push-pull based distributed algorithm to solve the dual problem \eqref{D}. Let us unfold the details.
	\subsection{Design of Distributed Algorithm}
	Consider a centralized optimization problem with equality and inequality constraints:
	\begin{equation}\label{CO}
		\min_{x\in\X}\tilde{f}(x)\quad {\rm s.t.} \quad  \tilde{h}(x)\in\K,
	\end{equation}
	where $\K\triangleq\R^p_-\times\{0_q\}$. It is easy to verify that the saddle point problem  associated with problem \eqref{CO} is defined by 
	\begin{equation}\label{CD}
		\max_{\lambda\in\K^\circ}\min_{x\in\X}\tilde{L}(x,\lambda)
	\end{equation}
	where $\tilde{L}(x,\lambda)=\tilde{f}(x)+\lambda^\top\tilde{h}(x)$. The schemes of the primal-dual projected gradient method \cite{Arr1960} for the  maximin problem are expressed as follows:
	\begin{align*}
		x^{k+1}&=P_{\X}\left(x^k-\alpha_k\nabla_x\tilde{L}(x^k,\lambda^k)\right),\\
		\lambda^{k+1}&=P_{\K^\circ}\left(\lambda^k+\alpha_k\nabla_\lambda\tilde{L}(x^k,\lambda^k)\right).
	\end{align*}
	Review the first-order optimality conditions (\cite{Roc1970}, Theorem 36.6) of problem \eqref{CO}, which indicates that point $x\in\X$ is a optimal solution to problem \eqref{CO} if and only if there is a  pair of points $(x,\lambda)$ such that 
	\begin{align*}
		x&=P_{\X}\left(x-\alpha\nabla_x\tilde{L}(x,\lambda)\right),\\
		\lambda&=P_{\K^\circ}\left(\lambda+\alpha\nabla_\lambda\tilde{L}(x,\lambda)\right)
	\end{align*}
	for any scalar $\alpha>0$. Hence the primal-dual projected gradient method also can be viewed as a fixed point method.\\
	\textbf{Push-Pull Method:}  Over directed time-varying multi-agent networks, consider an unconstrained distributed optimization problem: $$\min_{x\in\R^n}f(x)=\sum_{i=1}^mf_i(x)$$ 
	where $f_i$ is the local objective function of agent $i$. Paper \cite{Saad20} proposed the time-varying push-pull gradient method, where the specific updating formulas of variables are formulated as follows:
	\begin{align*}
		x_i^{k+1}&=\sum_{j=1}^mA(k)_{ij}x_j^k-\alpha y_i^k,\\
		y_i^{k+1}&=\sum_{j=1}^mB(k)_{ij}y_j^k+\nabla f_i(x_i^{k+1})-\nabla f_i(x_i^k),
	\end{align*}
	The row-stochastic matrix $A(k)$ is used to mix the estimates $x_j^k$ of the optimal point, and the column-stochastic matrix $B(k)$ is employed to track the average gradient of the global objective function. Recall the saddle point problem \eqref{D} corresponding to the distributed constrained optimization problem \eqref{P}: $$\max_{\lambda\in\K^\circ}\min_{x\in\X}\sum_{i=1}^m\cL_i(x,\lambda)$$
	Combining the primal-dual projected method with the time-varying push-pull gradient method, we proposed the following primal-dual distributed algorithm for problem \eqref{D}:
	\begin{align}
		x_i^{k+1}&=P_{\X}\left(\sum_{j=1}^ma_{ij}(k)x_j^k-\alpha_kz_i^k\right),\\
		\lambda_i^{k+1}&=P_Q\left(\sum_{j=1}^ma_{ij}(k)\lambda_j^k+\alpha_ky_i^k\right),\\
		z_i^{k+1}&=\sum_{j=1}^mb_{ij}(k)z_j^k+\nabla_x\cL_i(x_i^{k+1},\lambda_i^{k+1})\notag\\
		&\hspace{6em} -\nabla_x\cL_i(x_i^k,\lambda_i^k),\\
		y_i^{k+1}&=\sum_{j=1}^mb_{ij}(k)y_j^k+h_i\left(x_i^{k+1}\right)-h_i\left(x_i^k\right).
	\end{align}
	where $z_i^0=\nabla_x\cL_i(x_i^0,\lambda_i^0), y_i^0=h_i(x_i^0), \forall i\in V$. Similar to the push-pull gradient method, the sequences $\{z_i^k\}$ and $\{y_i^k\}$ are used to track the average gradients of the Lagrange function $\cL(x,\lambda)$ with respect to $x$ and $\lambda$, respectively.
	\begin{rem}
		For the updating formulas of the multipliers $\lambda_i^k$, we employ the projection operator $P_Q$ instead of $P_{\K^\circ}$. The reason is that we consider the saddle point problem 
		\begin{equation}\label{CD1}
			\max_{\lambda\in Q}\min_{x\in\X}\tilde{L}(x,\lambda)
		\end{equation}
		instead of problem \eqref{CD}. In fact, one can obtain that problem \eqref{CD} and \eqref{CD1} have the same saddle points, e.g., see Lemma 3.1 in \cite{Zhu12} .
	\end{rem}
	\begin{algorithm}[htbp]
		\renewcommand{\algorithmicrequire}{\textbf{Input:}}
		\renewcommand{\algorithmicensure}{\textbf{Output:}}
		\caption{Push-Pull Based Distributed Primal-Dual Algorithm.}
		\label{alg1}
		\begin{algorithmic}[1]
			\REQUIRE
			Let $n$ be the number of iteration steps of termination.
			Initial point $(x_i^0,\lambda_i^0)\in\X\times\K^\circ$, \\
			$y_i^0=h_i(x_i^0),\ z_i^0=\nabla f_i(x_i^0)+\nabla h_i(x_i^0)\lambda_i^0$,\\ 
			$i=1,2,\dots,m$. Stepsize $\alpha_k>0$. Set $k=0$.\\
			\ENSURE $x_i^k,\lambda_i^k,i=1,2,\dots,m$.
			\WHILE  {$k\le n$}
			\FOR {$i=1,2,\dots,m$}
			\STATE 	$v_i^k=\sum_{j=1}^ma_{ij}(k)x_j^k$,\\
			\vspace{0.2cm}
			\STATE  $u_i^k=\sum_{j=1}^ma_{ij}(k)\lambda_j^k$,\\
			\vspace{0.2cm}
			\STATE 	$x_i^{k+1}=P_\X(v_i^k-\alpha_kz_i^k)$,\\
			\vspace{0.2cm}
			\STATE 	$\lambda_i^{k+1}=P_Q(u_i^k+\alpha_ky_i^k)$,\\
			\vspace{0.2cm}
			\STATE  $ y_i^{k+1}=\sum_{j=1}^mb_{ij}(k)y_j^k$ \\ 
			\vspace{0.2cm}
			$\hspace{3em} +h_i(x_i^{k+1})-h_i(x_i^k),$\\
			\vspace{0.2cm}
			\STATE 	$d_i^k=\nabla f_i(x_i^k)+\nabla h_i(x_i^k)\lambda_i^k$,\\
			\vspace{0.2cm}
			\STATE 	$d_i^{k+1}=\nabla f_i(x_i^{k+1})+\nabla h_i(x_i^{k+1})\lambda_i^{k+1}$,\\
			\vspace{0.2cm}
			\STATE 	$z_i^{k+1}=\sum_{j=1}^mb_{ij}(k)z_j^k+d_i^{k+1}-d_i^k$,
			\ENDFOR
			\STATE Update $k=k+1$,
			\ENDWHILE
		\end{algorithmic}
	\end{algorithm}
	\par In order to ensure that the sequence $\{(x_i^k,\lambda_i^k)\}$ generated by Algorithm \ref{alg1} could converge to a optimal solution to problem \eqref{D}, we require the stepsizes $\{\alpha_k\}$ to satisfy the following assumption.
	\begin{ass}\label{ass5}
		Suppose that the stepsizes $\{\alpha_k\}_{k=0}^\infty$ satisfy 
		\begin{equation}
			\sum_{k=0}^\infty\alpha_k=\infty,\quad \sum_{k=0}^\infty\alpha_k^2<\infty.
		\end{equation} 
	\end{ass}
	Under Assumption \ref{ass3}, the subset $\X$ is compact, according to the continuity and differentiablity of $f_i$ and $h_i$, the gradients of $f_i$ and $h_i$ are Lipschitz continuous over $\X$. Without loss of generality, there exists a constant $L$ such that
	\begin{equation}\label{Lc}
		\begin{aligned}
			\Vert f_i(x)\Vert\le L,\quad\Vert h_i(x)&\Vert\le L,\\ 
			\Vert\nabla f_i(x)-\nabla f_i(y)\Vert&\le L\Vert x-y\Vert,\\
			\Vert h_i(x)-h_i(y)\Vert&\le L\Vert x-y\Vert,\\
			\Vert\nabla h_i(x)-\nabla h_i(y)\nabla\Vert&\le L\Vert x-y\Vert,\\
			\forall x, y\in\X,& i\in V.
		\end{aligned}
	\end{equation}
	\section{Main Results}\label{sec:3}
	In this section, we will analyze the convergence and iteration complexity of \Cref{alg1}. Here we first show the obtained results in this paper as follows:
	\begin{thm}\label{thm:1}
		Under Assumptions \ref{ass1} - \ref{ass5}, for any $i\in\{1,2,\dots,m\}$, the sequence $\{(x_i^k,\lambda_i^k)\}$ generated by \Cref{alg1} converges to some saddle point $(x^*,\lambda^*)$ of the Lagrange function $\cL(x,\lambda)$, i.e., 
		\begin{equation}
			\begin{aligned}
				\lim_{k\to\infty}\Vert x_i^k-x^*\Vert&=0,\\
				\lim_{k\to\infty}\Vert\lambda_i^k-\lambda^*\Vert&=0, \forall i=1,2,\dots,m.
			\end{aligned}
		\end{equation}
	\end{thm}
	For the statement of next theorem, we first define the time-average estimate $\tilde{x}_s^n$ as follows:
	\begin{equation}\label{time-a-e}
		\tilde{x}_s^n\triangleq\frac{\sum_{k=s}^n\alpha_k\bar{x}^k}{\sum_{k=s}^n\alpha_k},\quad\forall \ 0\le s\le n.
	\end{equation}
	where $\bar{x}^k$ and $\bar{\lambda}^k$ are defined by equation \eqref{bar}. Then we have the following iteration complexity result for Algorithm \ref{alg1}.
	\begin{thm}\label{thm:2}
		For Algorithm \ref{alg1}, if the stepsizes $\{\alpha_k\}_{k=0}^\infty$ meet Assumption \ref{ass5}, then there exists a constant $M_1>0$ such that
		\begin{equation}
			f(\tilde{x}_s^n)-f(x^*)\le\frac{M_1}{\sum_{k=s}^n\alpha_k}, \forall \ 0\le s\le n.
		\end{equation}
		where $x^*$ is one of the optimal solutions to problem \eqref{P}.
	\end{thm}
	\begin{rem}
		The stepsizes can be chosen as $\alpha_k=\frac{c}{k^{1/2+\beta}}, \beta\in(0,1/2], c>0$, clearly such a stepsize selection satisfies Assumotion \ref{ass5}. For this stepsize selection, we can use Theorem \ref{thm:2} obtain the iteration complexty. Now we give the details. Let $\lfloor\alpha\rfloor$ denote the integer part of a positive real number $\alpha$, then it holds that
		\begin{align*}
			\sum_{k=\lfloor n/2\rfloor}^n\frac{1}{k^{1/2+\beta}}&\ge\int_{\lfloor n/2\rfloor}^{n+1}\frac{dt}{t^{1/2+\beta}}\\
			=\frac{2}{1-2\beta}&\left((n+1)^{\frac{1}{2}-\beta}-\lfloor n/2\rfloor^{\frac{1}{2}-\beta}\right).
		\end{align*}
		Let $\phi(x)=x^{1/2-\beta},\beta\in(0,1/2]$, it's easy to verify that $\phi(x)$ is a concave function, then we have $$\phi(n+1)-\phi(\lfloor n/2\rfloor)\ge\phi^{'}(n+1)(n+1-\lfloor n/2\rfloor),$$ i.e.,
		\begin{align*}
			&(n+1)^{\frac{1}{2}-\beta}-\lfloor n/2\rfloor^{\frac{1}{2}-\beta}\\
			\ge&\left(\frac{1}{2}-\beta\right)\frac{(n+1)-\lfloor n/2\rfloor}
			{(n+1)^{1/2+\beta}}\\
			\ge&\left(\frac{1}{2}-\beta\right)\frac{(n+1)-(n+1)/2}{(n+1)^{1/2+\beta}}\\
			=&\left(\frac{1}{4}-\frac{\beta}{2}\right)(n+1)^{\frac{1}{2}-\beta}
		\end{align*}
		From Theorem \ref{thm:2}, it yields that 
		\begin{equation}
			\begin{aligned}
				f\left(\tilde{x}_{\lfloor n/2\rfloor}^n\right)-f\left(x^*\right)
				&\le\frac{M_1}{\sum_{k=\lfloor n/2\rfloor}^n\frac{c}{k^{1/2+\beta}}}\\
				\vspace{0.2cm}
				&\le \frac{2M_1}{c(n+1)^{1/2-\beta}}.
			\end{aligned}
		\end{equation}
		In this case, we can obtain the iteration complexity of Algorithm \ref{alg1} as $\mathcal{O}\left(\frac{1}
		{(n+1)^{1/2-\beta}}\right),\beta\in(0,1/2)$.\qed
	\end{rem}
	
	\subsection{Convergence Analysis}
	For the purpose of simplicity, we make the following notation protocols,
	\begin{align*}
		\omega_i^k&=\left(x_i^k,\lambda_i^k\right),
		\eta_i^k=\left(z_i^k,-y_i^k\right),\\
		\nabla \cL_i(\omega_i^k)&=\left(d_i^k,-h_i(x_i^k)\right),
		\hat{v}_i^k=\left(v_i^k,u_i^k\right),\\
		\omega^k&=\left(\omega_1^k,\omega_2^k,\cdots,\omega_m^k\right)^\top\in\R^{2m\times(n+p+q)},\\
		\eta^k&=\left(\eta_1^k,\eta_2^k,\cdots,\eta_m^k\right)^\top\in\R^{2m\times(n+p+q)},\\
		\nabla \mL(\omega^k)&=\left(\nabla\cL_1(\omega_1^k),\cdots,\nabla\cL_m(\omega_m^k)\right)^\top,\\
		\quad e^k&=\left(e_1^k,e_2^k,\cdots,e_m^k\right)^\top,
	\end{align*}
	\begin{align*}
		\A_k&=\left(
		\begin{array}{cc}
			A(k) & 0_{m\times m}\\
			0_{m\times m} &A(k) 
		\end{array}
		\right)\in\R^{2m\times 2m},\\
		\quad \B_k&=\left(
		\begin{array}{cc}
			B(k) &0_{m\times m}\\
			0_{m\times m} & B(k)
		\end{array}
		\right)\in\R^{2m\times 2m}.
	\end{align*}
	where $0_{m\times m}$ is the zeros matrix in $\R^{m\times m}$, and  $$e_i^k=\omega_i^{k+1}-(\hat{v}_i^k-\alpha_k\eta_i^k).$$
	Let 
	\begin{align*}
		\begin{matrix}
			m\text{\ pairs} \\ \Omega= \overbrace{\X\times Q\times\X\times Q\times\dots\times\X\times Q},
		\end{matrix}
	\end{align*}
	then we can rewrite the updating formulas of variables in Algorithm \ref{alg1} into the following compact matrix form:
	\begin{align}
		\omega^{k+1}&=P_\Omega\left(\A_k\omega^k-\alpha_k\eta^k\right),\label{eq:1}\\
		\eta^{k+1}&=\B_k\eta^k+\nabla\mL(\omega^{k+1})-\nabla\mL(\omega^k),\label{eq:2}
	\end{align}
	where the matrices $\{\A_k\}$ and $\{\B_k\}$ are row-stochastic and column-stochastic, respectively. To proceed with the analysis, we make a state transformation: $s^k=V_k^{-1}\eta^k$, where $V_k=\mathrm{diag}(v_k)$ and $v_k$ are defined by  
	\begin{equation}\label{eq:3}
		v_{k+1}=\B_kv_k,\quad v_0=1_{2m}.
	\end{equation}
	Then equations \eqref{eq:1} and \eqref{eq:2} are equivalent to 
	\begin{equation}\label{eq:4}
		\omega^{k+1}=P_\Omega\left(\A_k\omega^k-\alpha_kV_ks^k\right)
	\end{equation}
	and
	\begin{equation}\label{eq:5}
		s^{k+1}=\rc_ks^k+V_{k+1}^{-1}\left(\nabla\mL(\omega^{k+1})-\nabla\mL(\omega^k)\right),
	\end{equation}
	respectively, 
	where $\rc_k=V_{k+1}^{-1}\B_kV_k$. It can be verified that $\{\rc_k\}$ is a sequence of row-stochastic matrices for which the absolute probability sequence is $\{\frac{1}{2m}v_k\}$.
	Let us denote
	\begin{align*} \
		\Phi_\A(k,s)&=\A_k\times\A_{k-1}\times\cdots\times\A_s,\\
		\Phi_\rc(k,s)&=\rc_k\times\rc_{k-1}\times\cdots\times\rc_s, \forall k\ge s.
	\end{align*}
	The asymptotical states of matrix sequences $\{\Phi_\A(k,s)\}$ and $\{\Phi_\rc(k,s)\}$ are revealed by the following lemma.
	\begin{lem}[\cite{Nedic09a}, Lemma 4]\label{lem:2}
		Under Assumption \ref{ass1} and Assumption \ref{ass2}, then we have:
		\begin{enumerate}[\hspace{1em}(1)]
			\item The limit $\bar{\Phi}_\A(s)=\lim_{k\to\infty}\Phi_A(k,s)$ exists for each $s$.
			\item $\bar{\Phi}_\A(s)=1_{2m}\mu_s^\top$, where $\mu_s\in\R^{2m}$ is a stochastic vector for each $s$.
			\item For any $i$, the entrices $\Phi_\A(k,s)_{ij}, j=1,\dots,2m$, converge to the same limit $(\mu_s)_i$ as $k\to\infty$ with a geometric rate, i.e., for each $i$ and $s\ge0$,
			\begin{equation*}
				\left\vert\Phi_\A(k,s)_{ij}-\left(\mu_s\right)_i\right\vert\le C\beta^{k-s},\ \forall k\ge s,
			\end{equation*} 
			where 
			\begin{align*}
				C&=2\frac{1+2\eta^{-B_0}}{1-2\eta^{B_0}}\left(1-\eta\right)^{1/{B_0}},\\
				\beta&=(1-\eta^{B_0})\in(0,1), 
			\end{align*}
			$\eta$ is the lower bound of Assumption \ref{ass2} (3), $B_0=(m-1)B$, ingeter $B$ is defined by Assumption \ref{ass1}.
		\end{enumerate}
	\end{lem}
	The above lemma is also true if $\Phi_\A(k,s)$ is replaced by $\Phi_\rc(k,s)$, due to that the matrices $\{\Phi_\rc(k,s)\}$ are the same row-stochastic as $\{\Phi_\A(k,s)\}$. To prove Theorem \ref{thm:1} and Theorem \ref{thm:2}, we need some preparations, in which the following three lemmas are indispenable. Let us give the details.
	\begin{lem}[\cite{Saad20},Corollary 1]\label{lem:3}
		Under the assumptions of Lemma \ref{lem:2}, the sequence $\{\phi_k\}$ is an absolute probability sequence for the matrix sequence $\{\A_k\}$, where
		\begin{equation}
			\begin{aligned}
				\phi_k^\top&=\mu_s^\top,\quad k=sB,\\
				\phi_k^\top&=\mu_{s+1}^\top\A_{(s+1)B-1}\cdots\A_k,\\
				k&\in\left(sB,(s+1)B\right),\\
				\text{for} \ \ s&=0,1,2,\cdots.
			\end{aligned}
		\end{equation}
	\end{lem}
	\begin{lem}[\cite{Polyak20}, Lemma 11,Chapter 2.2]\label{lem:4}
		Let $\{b_k\}, \{c_k\}, \{d_k\}$ be non-negative sequences. Suppose that $\sum_{k=0}^\infty c_k<\infty$ and
		\begin{equation*}
			b_{k+1}\le b_k-d_k+c_k,\quad \forall k\ge1
		\end{equation*}
		then the sequence $\{b_k\}$ converges and $\sum_{k=0}^\infty d_k<\infty$.
	\end{lem}
	\begin{lem}[\cite{Nedic10}, Lemma 7]\label{lem:5}
		Let $\beta\in(0,1)$, and $\{\gamma_k\}$ be a positive scalar sequence.
		\begin{enumerate}[\hspace{1em}(1)]
			\item If $\lim_{k\to\infty}\gamma_k=0$, then
			$$\lim_{k\to\infty}\sum_{l=0}^k\beta^{k-l}\gamma_l=0,$$	
			\item In addition, if $\sum_{k=0}^\infty\gamma_k<\infty$, then
			$$\sum_{k=0}^\infty\sum_{l=0}^k\beta^{k-l}\gamma_l<\infty.$$
		\end{enumerate}
	\end{lem}
	From equation \eqref{eq:2}, we have
	\begin{align*}
		1_{2m}^\top\eta^{k+1}&=1_{2m}^\top\B_k\eta_k+1_{2m}^\top\left(\nabla\mL(\omega^{k+1})-\nabla\mL(\omega^k)\right)\\
		&=1_{2m}^\top\eta^k+1_{2m}^\top\left(\nabla\mL(\omega^{k+1})-\nabla\mL(\omega^k)\right),
	\end{align*}  
	with the given initial vectors $\omega_i^0=(x_i^0,\lambda_i^0), i=1,\dots,m$, it can be verified that $$1_{2m}^\top\eta^k=1_{2m}^\top\nabla\mL(\omega^k),\quad\forall k\ge0.$$
	Let
	\begin{equation}\label{bar}
		\begin{aligned}
			\bar{\omega}^k=\left(\bar{x}^k,\bar{\lambda}^k\right)
			=\phi_k^\top\omega^k,\quad \bar{s}^k=\frac{1}{2m}v_k^\top s^k,
		\end{aligned}
	\end{equation}
	where the sequence $\{\phi_k\}$ is defined by Lemma \ref{lem:3}.
	Then  we obtain 
	\begin{align*}
		\bar{s}^k&=\frac{1}{2m}v_k^\top s^k=\frac{1}{2m}v_k^\top V_k^{-1}\eta^k\\
		&=\frac{1}{2m}1_{2m}^\top\eta^k=\frac{1}{2m}1_{2m}^\top\nabla\mL(\omega^k),\quad \forall k\ge0,
	\end{align*}
	which means that $\bar{s}^k$ tracks the average gradient of the Lagrange function $\cL(x,\lambda)$.
	Since the subset $\X$ is compact convex, $\nabla_x\cL(x,\lambda)$ and $\nabla_\lambda\cL(x,\lambda)$ are bounded over $\X$, which implies the sequence $\bar{s}^k$ is bounded. In addition, the sequence $s^k$ is also bounded. We give the result in the following lemma. 
	\begin{lem}\label{lem:6}
		There is a constant $M>0$ such that $$\Vert s^k\Vert\le M,\quad \forall k\ge0.$$
	\end{lem} 
	For the  proof of this lemma, see Appendix \ref{A1}. Let us introduce another extremely important lemma for the convergence analysis in this article, which reveals the asymptotical relations between $\omega^k$ and $\bar{\omega}^k$, as well as $s^k$ and $\bar{s}^k$.
	\begin{lem}\label{lem:7}
		Under Assumption \ref{ass1} - \ref{ass5}, for the sequences generated by the Algorithm \ref{alg1}, it meets that
		\begin{align*}
			\lim_{k\to\infty}\Vert\omega^k-1_{2m}\bar{\omega}^k\Vert&=0,\\
			\lim_{k\to\infty}\Vert s^k-1_{2m}\bar{s}^k\Vert&=0,\\
			\sum_{k=0}^\infty\alpha_k\Vert\omega^k-1_{2m}\bar{\omega}^k\Vert&<\infty,\\
			\sum_{k=0}^\infty\alpha_k\Vert s^k-1_{2m}\bar{s}^k\Vert&<\infty.
		\end{align*}
	\end{lem}
	The proof of this lemma is given in Appendix \ref{A2}. Now we start to develop the proof of Theorem \ref{thm:1} and Theorem \ref{thm:2}.
	\par \textbf{Proof of Theorem \ref{thm:1}:} 
	For any $\omega^*=(x^*,\lambda^*)^\top$ $\in(\X\times\K^\circ)$, we define
	\begin{align*}
		I_1^k&=2\alpha_k\left\langle 1_{2m}(\nabla\mL(\omega^k)-\nabla \mL(1_{2m}\bar{\omega}^k)),\bar{\omega}^k-\omega^*\right\rangle,\\
		I_2^k&=2\alpha_k\left\langle V_ks^k-1_{2m}\bar{s}^k,\A_k\omega^k-1_{2m}\omega^*\right\rangle,\\
		I_3^k&=2\alpha_k\left\langle 1_{2m}(\bar{\omega}^k-\omega^*),\A_k\omega^k-1_{2m}\bar{\omega}^k\right\rangle.
	\end{align*}
	Then we split the proof into the following three steps. 
	\begin{itemize}
		\item \textbf{Step 1:} Verify the inequality
		\begin{equation}\label{ieqimpt}
			\begin{aligned}
				&\Vert\omega^{k+1}-1_{2m}\omega^*\Vert^2\\
				\le&\Vert\omega^k-1_{2m}\omega^*\Vert^2-
				2\alpha_k\left(\cL(\bar{x}^k,\lambda^*)-\cL(x^*,\lambda^*)\right)\\
				&\quad-2\alpha_k\left(\cL(x^*,\lambda^*)-\cL(x^*,\bar{\lambda}^k)\right)\\
				&\quad+I_1^k+I_2^k+I_3^k+4m^2M^2\alpha_k^2.
			\end{aligned}
		\end{equation}
		From equation \eqref{eq:4} and the Lipschitz continuity of projection operator $P_\Omega$, we have
		\begin{align*}
			&\Vert\omega^{k+1}-1_{2m}\omega^*\Vert^2\\
			\le&\Vert\A_k\omega^k-\alpha_kV_ks^k-1_{2m}\omega^*\Vert^2\\
			=&\Vert\A_k\omega^k-1_{2m}\omega^*\Vert^2+\alpha_k^2\Vert V_ks^k\Vert^2\\
			&\quad-2\alpha_k\left\langle V_ks^k,\A_k\omega^k-1_{2m}\omega^*\right\rangle
		\end{align*}
		Due to that the sequence $\{\frac{1}{2m}v_k\}$ are stochastic vectors, and $V_k=\mathrm{diag}(v_k)$, it follows that $$\Vert V_ks^k\Vert\le\Vert V_k\Vert\Vert s^k\Vert\le 2m\Vert s^k\Vert\le 2mM,$$ the last inequality comes from Lemma \ref{lem:6}. Since $\A_k$ is row-stochastic matrix, it holds that 
		$\Vert\A_k\omega^k-1_{2m}\omega^*\Vert^2\le\Vert\omega^k-1_{2m}\omega^*\Vert^2$.
		For the term of inner product, one has
		\begin{align*}
			&2\alpha_k\left\langle V_ks^k,\A_k\omega^k-1_{2m}\omega^*\right\rangle\\
			=&2\alpha_k\left\langle 1_{2m}\bar{s}^k,1_{2m}\bar{\omega}^k-1_{2m}\omega^*\right\rangle\\
			& +2\alpha_k\left\langle V_ks^k-1_{2m}\bar{s}^k,\A_k\omega^k-1_{2m}\omega^*\right\rangle\\
			& +\left\langle1_{2m}\bar{\omega}^k-\omega^*,\A_k\omega^k-1_{2m}\bar{\omega}^k\right\rangle\\
			=&\left\langle1_{2m}^\top\nabla\mL(1_{2m}\bar{\omega}^k),\bar{\omega}^k-\omega^*\right\rangle+I_1^k+I_2^k+I_3^k
		\end{align*}
		Note that the Lagrange function $\cL(x,\lambda)$ is convex-concave with respect to $x$ and $\lambda$, according to the subdifferential inequality, we obtain 
		\begin{align*}
			&\left\langle1_{2m}^\top\nabla\mL(1_{2m}\bar{\omega}^k),\bar{\omega}^k-\omega^*\right\rangle\\
			\ge&\left(\cL(\bar{x}^k,\lambda^*)-\cL(x^*,\bar{\lambda}^k)\right).
		\end{align*}
		Combining the above three inequalities, one can get inequality \eqref{ieqimpt}.
		\item \textbf{Step 2:} Verify that
		\begin{align*}
			&\sum_{k=0}^\infty\vert I_1^k\vert<\infty,\\
			&\sum_{k=0}^\infty\vert I_2^k\vert<\infty,\\
			&\sum_{k=0}^\infty\vert I_3^k\vert<\infty.
		\end{align*}
		According to the Cauchy inequality and the Lipschitz continuity of $\nabla\cL$, the following inequalities are obvious
		\begin{align*}
			\vert I_1^k\vert&\le2L_1\alpha_k\Vert\omega^k-1_{2m}\bar{\omega}^k\Vert,\\
			\vert I_2^k\vert&\le2L_2\alpha_k\Vert s^k-1_{2m}\bar{s}^k\Vert,\\
			\vert I_3^k\vert&\le2L_3\alpha_k\Vert \omega^k-1_{2m}\bar{\omega}^k\Vert
		\end{align*}
		for some positive  constants $L_1, L_2$ and $L_3$. Then the convergence of these three series is from Lemma \ref{lem:7}.
		\item \textbf{Step 3:} Prove the convergence. In inequality \eqref{ieqimpt}, we take $\omega^*=(x^*;\lambda^*)$ as a saddle point of function $\cL(x,\lambda)$. By the definition of saddle point, it yields that
		\begin{align*}
			\cL(\bar{x}^k,\lambda^*)-\cL(x^*,\lambda^*)&\ge0,\\
			\cL(x^*,\lambda^*)-\cL(x^*,\bar{\lambda}^k)&\ge0.
		\end{align*}
		Therefore, it follows from Lemma \ref{lem:4} that the limit $\lim_{k\to\infty}\Vert\omega^k-1_{2m}\omega^*\Vert$ exists, and
		\begin{align*}
			\sum_{k=0}^\infty\alpha_k\left(\cL(\bar{x}^k,\lambda^*)-\cL(x^*,\lambda^*)\right)&<\infty,\\
			\sum_{k=0}^\infty\alpha_k\left(\cL(x^*,\lambda^*)-\cL(x^*,\bar{\lambda}^k)\right)&<\infty.
		\end{align*}
		As the stepsizes satisfy that $\sum_{k=0}^\infty\alpha_k=\infty$, we obtain
		\begin{align*}
			\liminf_{k\to\infty}\cL(\bar{x}^k,\lambda^*)&=\cL(x^*,\lambda^*),\\
			\limsup_{k\to\infty}\cL(x^*,\bar{\lambda}^k)&=\cL(x^*,\lambda^*).
		\end{align*}
		Since the Lagrange function $\cL(x,\lambda)$ is continuous differentiable and convex-convex with respect to $x$ and $\lambda$, it holds that
		\begin{equation*}
			\lim_{k\to\infty}\Vert\bar{x}^k-x^*\Vert=0,\quad\lim_{k\to\infty}\Vert\bar{\lambda}^k-\lambda^*\Vert=0.
		\end{equation*}
		From Lemma \ref{lem:7}, it follows that $$\lim_{k\to\infty}\Vert\omega^k-1_{2m}\omega^*\Vert=0,$$ 
		which is exactly the conclusion to be proved.\qed
	\end{itemize}
	
	\par \textbf{Proof of Theorem \ref{thm:2}:} For any fixed $\lambda\in\K^\circ$, the function $\cL(\cdot,\lambda)$ is convex, by the definition of $\tilde{x}_s^n$, we have 
	\begin{align*}
		&f(\tilde{x}_s^n)-f(x^*)\\
		\le&\cL(\tilde{x}_s^n,\lambda^*)-\cL(x^*,\lambda^*)\\
		\le&\left(\sum_{k=s}^n\alpha_k\right)^{-1}\left(\sum_{k=s}^n\alpha_k
		\left(\cL(\bar{x}^k,\lambda^*)-\cL(x^*,\lambda^*)\right)\right).
	\end{align*} 
	From the proof of Theorem \ref{thm:1}, we obtained that $$0<\sum_{k=0}^\infty\alpha_k\left(\cL(\bar{x}^k,\lambda^*)-\cL(x^*,\lambda^*)\right)=M_1<\infty.$$
	Therefore, it is clear that
	\begin{equation*}
		f(\tilde{x}_s^n)-f(x^*)\le\frac{M_1}{\sum_{k=s}^n\alpha_k}.
	\end{equation*}\qed
	\section{Simulation}\label{sec:4}
\begin{figure}[htbp]
		\centering
		\includegraphics[width=0.5\textwidth]{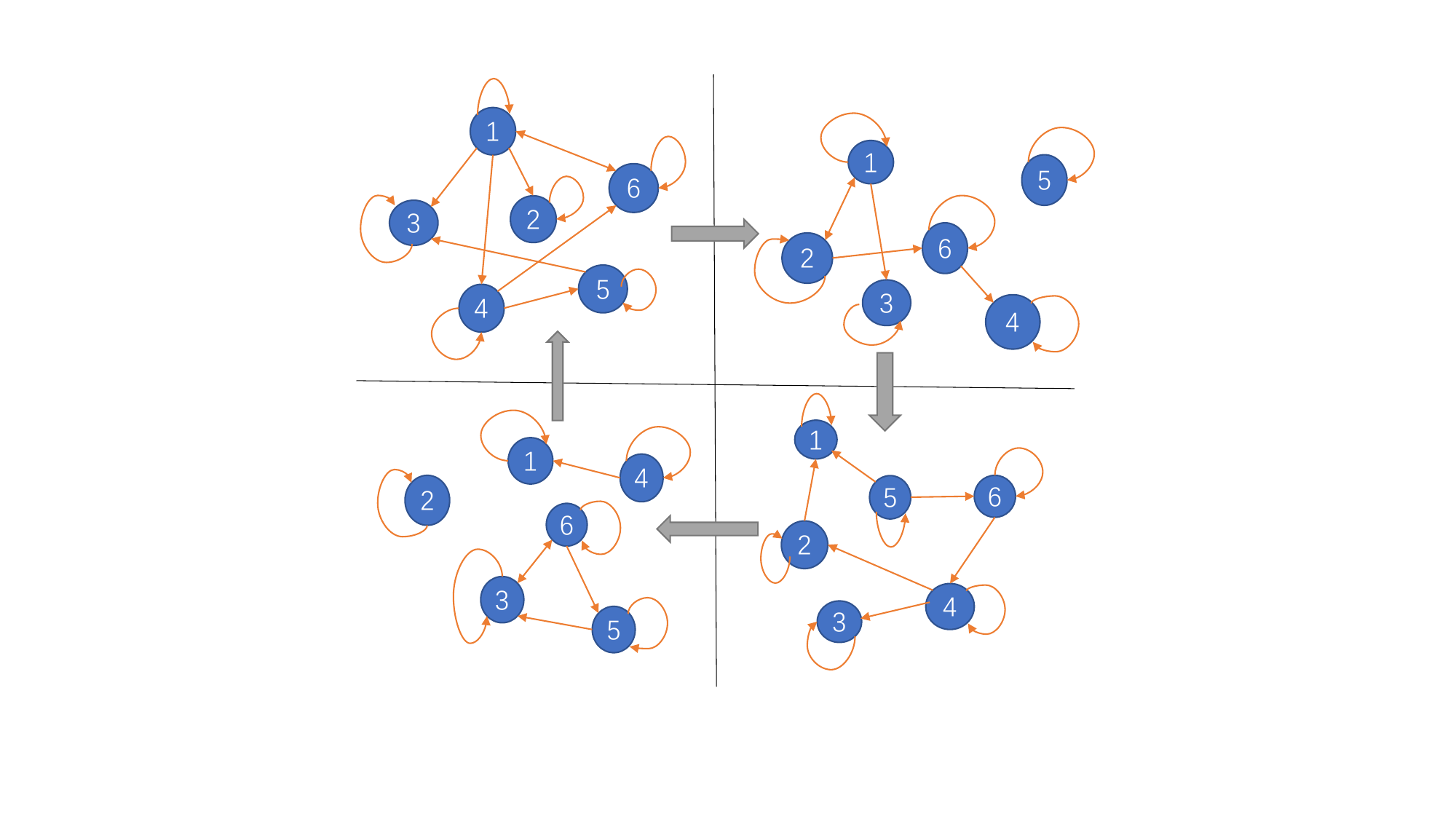}
		\caption{Periodic time-varying network topology, where the period is 4.}\label{Fig:1}
	\end{figure}

\begin{figure}[htbp]
		\centering
		\includegraphics[width=0.5\textwidth]{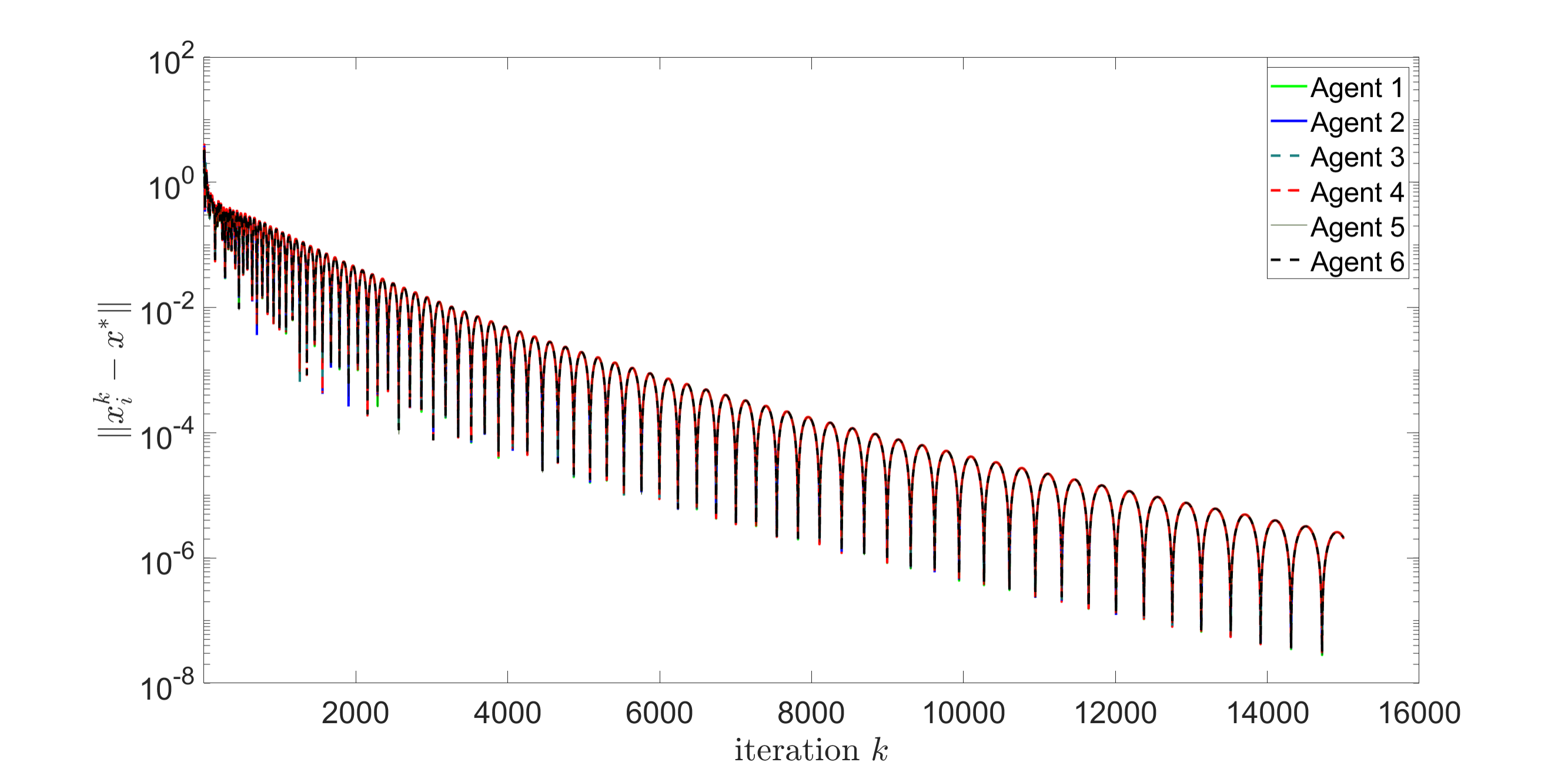}
		\vspace{-0.6cm}
		\caption{Evolution of the  distance between the decision variable $x_i^k$ of agent $i$ and the optimal solution $x^*$ to problem \eqref{P}, $i\in V$ .}\label{Fig:2}
	\end{figure}

\begin{figure}[htbp]
		\centering
		\includegraphics[width=0.5\textwidth]{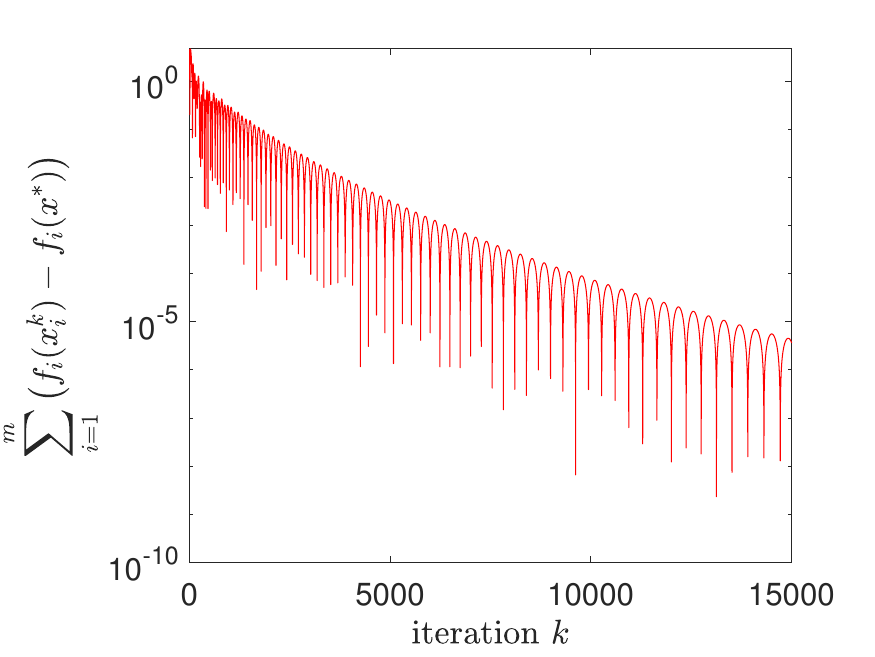}
		\vspace{-0.6cm}
		\caption{Evolution of the diffierence betweent the values of the global objective function at $x_i^k$ and $x^*$, where $x^*$ is a optimial solution to problem \eqref{P}.}\label{Fig:3}
	\end{figure}

\begin{figure}[htbp]
		\centering
		\includegraphics[width=0.5\textwidth]{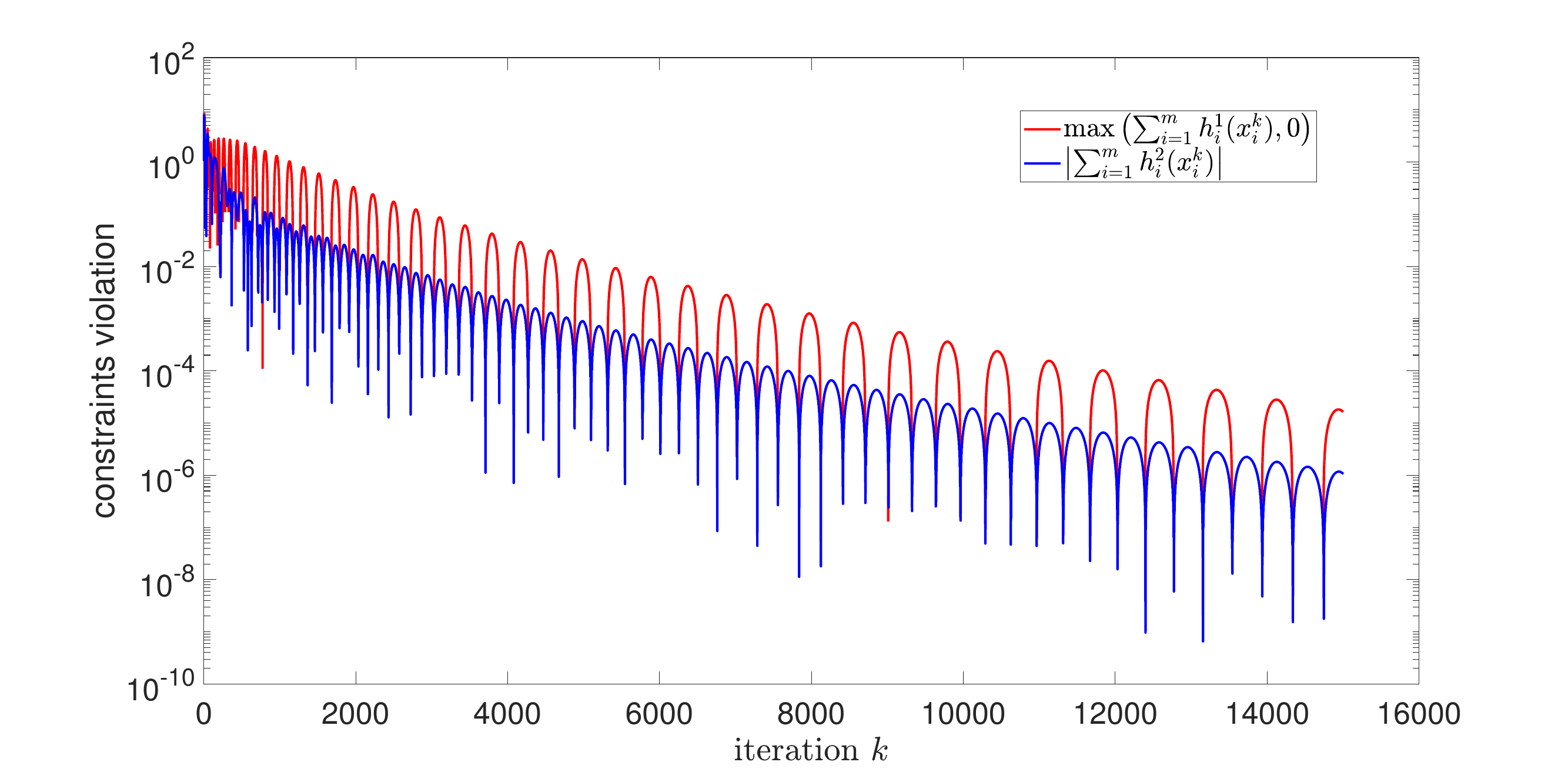}
		\caption{Evolution of constraints violation $\left[\max\left(\sum_{i=1}^mh_i^1(x),0\right), \left\vert\sum_{i=1}^mh_i^2(x)\right\vert\right]$ as a function of $x_j^k, j\in V$.}\label{Fig:4}
	\end{figure}

	In this section, we illustrate Algorithm \ref{alg1} by an example. Consider problem \eqref{P} with 
	\begin{align*}
		f_i(x)&=a_i^\top x+b_i+c_i\log(1+e^{d_i^\top x}),\\
		h_i^1(x)&=\alpha_i\left\Vert x\right\Vert^2+\beta_i,\quad p=1,\\
		h_i^2(x)&=\gamma_i^\top x+\delta_i, \quad q=1,\\
		\X&=[-3,3]\times [-3,3].
	\end{align*}
	This problem was also considered in \cite{Liang20}, while the local objective functions of agents in the example of \cite{Liang20} are strongly convex and L-smooth, here they are just convex and smooth. The inequality consraints decribe an elliptic region, and the equality constraints characterize an affine subspace. It is easy to verify that the Slater's condition holds in Assumption \ref{ass4}. We take $m=6$ and use directed time-varying graphs. The network topologies vary according to the periodic sequence of directed graphs as shown in \Cref{Fig:1} making the directed communication network 4-bounded strongly connected. Then we adopt a simple uniform weighting strategy to construct the row-stochastic and column-stochastic weight matrices $A(k)$ and $B(k)$ as follows:
	\begin{align*}
		a_{ij}(k)&=\left\{
		\begin{aligned}
			\frac{1}{\vert\N_i^{\rm{in}}(k)\vert}&,\quad (i,j)\in E(k),\\
			0\hspace{1.2em}&,\quad (i,j)\notin E(k),
		\end{aligned}
		\right.
		\\
		b_{ij}(k)&=\left\{
		\begin{aligned}
			\frac{1}{\vert\N_j^{\rm{out}}(k)\vert}&,\quad (i,j)\in E(k),\\
			0\hspace{1.2em}&,\quad (i,j)\notin E(k),
		\end{aligned}
		\right.
	\end{align*}     
	where $\vert\N_i^{\rm{in}}(k)\vert$ and $\vert\N_j^{\rm{out}}(k)\vert$ are the in-degree of agent $i$ and out-degree of agent $j$ at time $k$, respectively. In the simulation, we select the stepsizes as $\alpha_k=\frac{2}{(k+1)^{0.6}}$. Then our procedures put out the following images \Cref{Fig:2}, \Cref{Fig:3} and \Cref{Fig:4}.
	\par \Cref{Fig:2} indicates that the decision variables of all agents will reach a consensus and converge to a optimal solution to problem \eqref{P}. \Cref{Fig:3} shows that our algorithm does push all agents to adjust their decisions towards the goal of minimizing the global objective function. \Cref{Fig:4} suggests that at the beginning of the procedure, the decision variables $x_j^k$ of all agents do not satisfy the coupled constraints. After adjustment, the iterates render the coupled constraints and the optimum. 
	
	\section{Conclusion}\label{sec:5}
	In this article, we considered a distributed convex optimization problem with coupled equality and inequality constraints and a common constraint set in directed unbalanced and time-varying multi-agent networks. Based on the push-pull method for unconstrained case, we proposed a distributed primal-dual algorithm to solve the involoved problem, established the asymptotical convergence, and analyzed the iteration complexity of  our algorithm. In simulation, we showed that our algorithm is effictive for the problem considered in this papaer. However, there are still some respects that need to improve. Obviously, one of them is the stepsize, we adopted a diminishing stepsize, which leads to the slow covergence. Therefore, we expect to improve the diminishing stepsize to fixed stepsize in the future studies.
	
	\appendix
	\section{Proof of Lemma \ref{lem:6}}\label{A1}
	For the simplicity, denote 
	\begin{equation*}
		\delta_k=V_{k+1}^{-1}\left(\nabla\mL(\omega^{k+1})-\nabla\mL(\omega^k)\right),
	\end{equation*}
	then equation \eqref{eq:5} is equivalently  written as:
	\begin{equation*}
		s^{k+1}=\rc_ks^k+\delta_k.
	\end{equation*}
	Moreover, 
	\begin{equation*}
		s^{k+1}=\Phi_\rc(k,1)s^1+\sum_{l=1}^{k-1}\Phi_\rc(k,l+1)\delta_l+\delta_k,
	\end{equation*}
	Since $\{\frac{1}{2m}v_k\}$ is an absolute probability sequence for the matrix sequence $\{\rc_k\}$, i.e., $$\frac{1}{2m}v_k^\top=\frac{1}{2m}v_{k+1}^\top\rc_k,$$
	Then it holds that
	\begin{align*}
		&1_{2m}\bar{s}^{k+1}=\frac{1}{2m}1_{2m}v_{k+1}^\top s^{k+1}\\
		&=\frac{1}{2m}1_{2m}v_{k+1}^\top\left(\Phi_\rc(k,1)s^1+\delta_k\right)\\
		&\quad+\frac{1}{2m}1_{2m}v_{k+1}^\top\sum_{l=1}^{k-1}\Phi_\rc(k,l+1)\delta_l\\
		&=\frac{1}{2m}1_{2m}v_1^\top s^1+\frac{1}{2m}\sum_{l=1}^{k-1}1_{2m}v_l^\top\delta_l\\
		&\quad+\frac{1}{2m}1_{2m}v_{k+1}^\top\delta_k
	\end{align*}
	From the triangle inequality, we have 
	\begin{equation}\label{eq:6}
		\begin{aligned}
			&\Vert s^{k+1}-1_{2m}\bar{s}^{k+1}\Vert\\
			\le&\left\Vert\left(\Phi_\rc(k,1)-\frac{1}{2m}1_{2m}v_1^\top\right)s^1\right\Vert\\
			+&\left\Vert\left(I_{2m}-\frac{1}{2m}1_{2m}v_{k+1}^\top\right)\delta_k\right\Vert\\
			+&\sum_{l=1}^{k-1}\left\Vert\left(\Phi_\rc(k,l+1)-\frac{1}{2m}1_{2m}v_l^\top\right)\delta_l\right\Vert
		\end{aligned}
	\end{equation}
	For matrix sequence $\{\Phi_\rc(k,s)\}$, there holds the similar results in Lemma \ref{lem:4}, i.e., $$\left\Vert\Phi_\rc(k,s)-\frac{1}{2m}1_{2m}v_s^\top\right\Vert\le C_1\rho^{k-s}, \forall k\ge s,$$ where $\rho\in(0,1)$ and $C_1$ is some suitable constant. Due to the compactness of the subset $\X$, it is clear that $\delta_k$ is bounded. Assume that the upper bound of $\Vert\delta_k\Vert$ is a constant $\widetilde{C}$. Then it follows that
	\begin{align*}
		&\Vert s^{k+1}-1_{2m}\bar{s}^{k+1}\Vert\\
		\le&\rho^{k-1}\Vert s^1\Vert+\widetilde{C}\sum_{l=1}^{k-1}\rho^{k-1-l}+\widetilde{C}\\
		\le& \Vert s^1\Vert +\frac{\widetilde{C}}{1-\rho}+\widetilde{C}=:M,\quad \forall k\ge0.
	\end{align*}
	\qed
	\section{Proof of Lemma \ref{lem:7}}\label{A2}
	To achieve the goal of proving Lemma \ref{lem:7}, we need to verify the following inequality 
	\begin{align*}
		&\Vert\omega^{k+1}-1_{2m}\bar{\omega}^{k+1}\Vert\\ \le&C_2\beta^{k-1}+C_3\sum_{l=1}^{k-1}\beta^{k-1-l}\alpha_l+C_4\alpha_k.
	\end{align*}
	Then combining the conditions on stepsizes $\{\alpha_k\}$ in Assumption \ref{ass5} and Lemma \ref{lem:5}, it is easy to find the requied results.\\
	\textbf{Part (a):} From equation \eqref{eq:4}, $\omega^{k+1}=\A_k\omega^k-\alpha_kV_ks^k+e^k$. Similar to inequality \eqref{eq:6}, for the sequence $\{\omega^k\}$, we have 
	\begin{equation}\label{eq:7}
		\begin{aligned}
			&\Vert\omega^{k+1}-1_{2m}\bar{\omega}^{k+1}\Vert\\
			\le&\left\Vert\left(\Phi_\A(k,1)-1_{2m}\phi_1^\top\right)\omega^1\right\Vert\\
			+&\left\Vert\left(I_{2m}-1_{2m}\phi_k^\top\right)\left(e^k-\alpha_kV_ks^k\right)\right\Vert\\
			+&\sum_{l=1}^{k-1}\left\Vert\left(\Phi_\A(k,l+1)-1_{2m}\phi_l^\top\right)\left(e^l-\alpha_lV_ls^l\right)\right\Vert
		\end{aligned}
	\end{equation}
	Since $\omega^{k+1}$ is the projected point of $\A_k\omega^k-\alpha_kV_ks^k$ onto $\Omega$. Then for any $\omega\in\Omega$, we have 
	\begin{align*}
		&-\left\langle e^k, \omega-\omega^{k+1}\right\rangle\\
		=&\left\langle\A_k\omega^k-\alpha_kV_ks^k-\omega^{k+1},\omega-\omega^{k+1}\right\rangle\le0,
	\end{align*}
	which implies that
	\begin{align*}
		&\Vert\omega^{k+1}-\A_k\omega^k\Vert^2\\
		=&\Vert e^k-\alpha_kV_ks^k\Vert^2\\
		=&\Vert e^k\Vert^2-2\left\langle\alpha_kV_ks^k,e^k\right\rangle+\Vert\alpha_kV_ks^k\Vert^2\\
		=&\Vert e^k\Vert^2+2\left\langle\omega^{k+1}-\A_k\omega^k-e^k,e^k\right\rangle+\Vert\alpha_kV_ks^k\Vert^2\\
		\le&\alpha_k^2\Vert V_ks^k\Vert^2-\Vert e^k\Vert^2.
	\end{align*}
	From Lemma \ref{lem:6}, it holds that 
	\begin{equation*}
		\Vert e^k\Vert\le\alpha_k\Vert V_ks^k\Vert\le\alpha_k\Vert s^k\Vert\le M\alpha_k, \forall k\ge0.
	\end{equation*}
	Based on the results of Lemma \ref{lem:4} and Lemma \ref{lem:5}, we obtain an improved version of inequality \eqref{eq:7} as follows:
	\begin{equation}\label{eq:8}
		\begin{aligned}
			&\Vert\omega^{k+1}-1_{2m}\bar{\omega}^{k+1}\Vert\\
			\le&C\Vert\omega^1\Vert\beta^{k-1}+2M\alpha_k+
			2M\sum_{l=1}^{k-1}\beta^{k-1-l}\alpha_l
		\end{aligned}
	\end{equation}
	Under Assumption \ref{ass5} and Lemma \ref{lem:5} (1), it follows that 
	$$\lim_{k\to\infty}\Vert\omega^k-1_{2m}\bar{\omega}^k\Vert=0.$$
	\textbf{Part (b):} Multiplying inequality \eqref{eq:8} by $\alpha_{k+1}$ at the both sides, it yields that
	\begin{align*}
		&\alpha_{k+1}\Vert\omega^{k+1}-1_{2m}\bar{\omega}^{k+1}\Vert\\
		\le& C_2\beta^{k-1}\alpha_{k+1}+2C_3\sum_{l=1}^{k-1}
		\beta^{k-1-l}\alpha_l\alpha_{k+1}+C_4\alpha_k\alpha_{k+1}\\
		\le& C_2\left(\beta^{2(k-1)}+\alpha_{k+1}^2\right)+
		C_3\sum_{l=1}^{k-1}\beta^{k-1-l}\left(\alpha_l^2+\alpha_{k+1}^2\right)\\
		&\quad+C_4\left(\alpha_k^2+\alpha_{k+1}^2\right)\\
		\le& C_1\beta^{2(k-1)}+C_2\sum_{l=1}^{k-1}\beta^{k-1-l}\alpha_l^2+C_4\alpha_k^2+C_5\alpha_{k+1}^2
	\end{align*}
	where $C_2, C_3, C_4, C_5$ are some suitable constants. Once again, under the Assumption \ref{ass5} and Lemma \ref{lem:5} (2), it is concluded that 
	$$\sum_{k=0}^\infty\alpha_k\Vert\omega^k-1_{2m}\bar{\omega}^k\Vert<\infty.$$
	For the sequence $\{s^k\}$, the desired results in Lemma \ref{lem:7} can be got by a similar deduction to $\{\omega^k\}$, which is not repeated here.\qed
	\bibliographystyle{unsrt}
	\bibliography{ref_dpdstm}

\begin{thebibliography}{10}

\bibitem{Mahmoudi20}
Afsaneh Mahmoudi, Hossein~S. Ghadikolaei, and Carlo Fischione.
\newblock Cost-efficient distributed optimization in machine learning over
  wireless networks.
\newblock In {\em ICC 2020 - 2020 IEEE International Conference on
  Communications (ICC)}, pages 1--7, 2020.

\bibitem{Liao21b}
Shengbin Liao.
\newblock A fast distributed algorithm for coupled utility maximization problem
  with application for power control in wireless sensor networks.
\newblock {\em Journal of Communications and Networks}, 23(4):271--280, 2021.

\bibitem{Sasso20}
Francesco Sasso, Angelo Coluccia, and Giuseppe Notarstefano.
\newblock Interaction-based distributed learning in cyber-physical and social
  networks.
\newblock {\em IEEE Transactions on Automatic Control}, 65(1):223--236, 2020.

\bibitem{Sharma16}
Harpreet Sharma and Gagandeep Kaur.
\newblock Optimization and simulation of smart grid distributed generation: A
  case study of university campus.
\newblock In {\em 2016 IEEE Smart Energy Grid Engineering (SEGE)}, pages
  153--157, 2016.

\bibitem{Nowak03}
R.D. Nowak.
\newblock Distributed em algorithms for density estimation and clustering in
  sensor networks.
\newblock {\em IEEE Transactions on Signal Processing}, 51(8):2245--2253, 2003.

\bibitem{Feng13}
Yuwu Feng and Qinghua Huang.
\newblock Fully distributed acoustic source localization in wireless sensor
  network.
\newblock In {\em 2013 IEEE International Conference of IEEE Region 10 (TENCON
  2013)}, pages 1--4, 2013.

\bibitem{Liao21a}
Shengbin Liao.
\newblock A fast distributed algorithm for coupled utility maximization problem
  with application for power control in wireless sensor networks.
\newblock {\em Journal of Communications and Networks}, 23(4):271--280, 2021.

\bibitem{Shaoqun16}
Shaoqun Song, Ruipeng Guo, Feng Chen, and Wenying Huang.
\newblock A real-time power flow optimal control method for hybrid ac/dc power
  systems with vsc-hvdc.
\newblock In {\em 2016 International Conference on Smart Grid and Electrical
  Automation (ICSGEA)}, pages 26--30, 2016.

\bibitem{Nedic09a}
A.~Nedi\'c and A.~Ozdaglar.
\newblock Distributed subgradient methods for multi-agent optimization.
\newblock {\em IEEE Transactions on Automatic Control}, 54(1):48--61, 2009.

\bibitem{Shi15a}
W.~Shi, Q.~Ling, G.~Wu, and W.~T. Yin.
\newblock Extra: An exact first-order algorithm for decentralized consensus
  optimization.
\newblock {\em SIAM Journal on Optimization}, 25(2):944--966, 2015.

\bibitem{Nedic17a}
A.~Nedi\'c, A.~Olshevsky, and W.~Shi.
\newblock Achieving geometric convergence for distributed optimization over
  time-varying graphs.
\newblock {\em Siam Journal on Optimization}, 27(4):2597--2633, 2017.

\bibitem{Qu18}
Guannan Qu and Na~Li.
\newblock Harnessing smoothness to accelerate distributed optimization.
\newblock {\em IEEE Transactions on Control of Network Systems},
  5(3):1245--1260, 2018.

\bibitem{Xu15}
Jinming Xu, Shanying Zhu, Yeng~Chai Soh, and Lihua Xie.
\newblock Augmented distributed gradient methods for multi-agent optimization
  under uncoordinated constant stepsizes.
\newblock In {\em 2015 54th IEEE Conference on Decision and Control (CDC)},
  pages 2055--2060, 2015.

\bibitem{Xu18}
Jinming Xu, Shanying Zhu, Yeng~Chai Soh, and L.~H. Xie.
\newblock Convergence of asynchronous distributed gradient methods over
  stochastic networks.
\newblock {\em IEEE Transactions on Automatic Control}, 63(2):434--448, 2018.

\bibitem{Pu18}
S.~Pu, W.~Shi, J.~M. Xu, and A.~Nedi\'c.
\newblock Push–pull gradient methods for distributed optimization in
  networks.
\newblock {\em IEEE Transactions on Automatic Control}, 66(1):1--16, 2020.

\bibitem{Saad20}
Fakhteh Saadatniaki, Ran Xin, and Usman~A. Khan.
\newblock Decentralized optimization over time-varying directed graphs with row
  and column-stochastic matrices.
\newblock {\em IEEE Transactions on Automatic Control}, 65(11):4769--4780,
  2020.

\bibitem{Du18}
Wen Du, Lisha Yao, Di~Wu, Xinrong Li, Guodong Liu, and Tao Yang.
\newblock Accelerated distributed energy management for microgrids.
\newblock In {\em 2018 IEEE Power Energy Society General Meeting (PESGM)},
  pages 1--5, 2018.

\bibitem{Xin18}
Ran Xin and Usman~A. Khan.
\newblock A linear algorithm for optimization over directed graphs with
  geometric convergence.
\newblock {\em IEEE Control Systems Letters}, 2(3):315--320, 2018.

\bibitem{Xin19}
Ran Xin, Chenguang Xi, and Usman~A. Khan.
\newblock Frost—fast row-stochastic optimization with uncoordinated
  step-sizes.
\newblock {\em EURASIP Journal on Advances in Signal Processing}, 2019(1):1,
  2019.

\bibitem{Nedic10}
A.~Nedi\'c, A.~Ozdaglar, and P.~A. Parrilo.
\newblock Constrained consensus and optimization in multi-agent networks.
\newblock {\em IEEE Transactions on Automatic Control}, 55(4):922--938, 2010.

\bibitem{Lei16}
J.~L. Lei, H.~F. Chen, and H.~T. Fang.
\newblock Primal–dual algorithm for distributed constrained optimization.
\newblock {\em Systems and Control Letters}, 96:110--117, 2016.

\bibitem{Zhu12}
M.~H. Zhu and S.~Martínez.
\newblock On distributed convex optimization under inequality and equality
  constraints.
\newblock {\em IEEE Transactions on Automatic Control}, 57(1):151--164, 2011.

\bibitem{Carli20}
Raffaele Carli and Mariagrazia Dotoli.
\newblock Distributed alternating direction method of multipliers for linearly
  constrained optimization over a network.
\newblock {\em IEEE Control Systems Letters}, 4(1):247--252, 2020.

\bibitem{Falsone17}
Alessandro Falsone, Kostas Margellos, Simone Garatti, and Maria Prandini.
\newblock Dual decomposition for multi-agent distributed optimization with
  coupling constraints.
\newblock {\em Automatica}, 84:149--158, 2017.

\bibitem{Chang14}
Tsung-Hui Chang, Angelia Nedić, and Anna Scaglione.
\newblock Distributed constrained optimization by consensus-based primal-dual
  perturbation method.
\newblock {\em IEEE Transactions on Automatic Control}, 59(6):1524--1538, 2014.

\bibitem{Liang20}
Shu Liang, Le~Yi Wang, and George Yin.
\newblock Distributed smooth convex optimization with coupled constraints.
\newblock {\em IEEE Transactions on Automatic Control}, 65(1):347--353, 2020.

\bibitem{Yang19}
Qing Yang and Gang Chen.
\newblock Primal-dual subgradient algorithm for distributed constraint
  optimization over unbalanced digraphs.
\newblock {\em IEEE Access}, 7:85190--85202, 2019.

\bibitem{Zhang20}
Jiawei Zhang, Songyang Ge, Tsung-Hui Chang, and Zhi-Quan Luo.
\newblock A proximal dual consensus method for linearly coupled multi-agent
  non-convex optimization.
\newblock In {\em ICASSP 2020 - 2020 IEEE International Conference on
  Acoustics, Speech and Signal Processing (ICASSP)}, pages 5740--5744, 2020.

\bibitem{Nedic18}
Angelia Nedi\'c, Alex Olshevsky, and Wei Shi.
\newblock Improved convergence rates for distributed resource allocation.
\newblock In {\em 2018 IEEE Conference on Decision and Control (CDC)}, pages
  172--177, 2018.

\bibitem{Nedic09b}
Angelia Nedi\'c and Asuman Ozdaglar.
\newblock Approximate primal solutions and rate analysis for dual subgradient
  methods.
\newblock {\em SIAM Journal on Optimization}, 19(4):1757--1780, 2009.

\bibitem{Arr1960}
K.~J. Arrow, L.~Hurwicz, and H.~Uzawa.
\newblock {\em Studies in Linear and Nonlinear Programming}.
\newblock Stanford University Press, Palo Alto, 1958.

\bibitem{Roc1970}
R.~T. Rockafellar.
\newblock {\em Convex Analysis}.
\newblock Princeton Mathematical Series, 1970.

\bibitem{Polyak20}
Boris Polyak.
\newblock {\em Introduction to Optimization}.
\newblock Optimization Software, 07 2020.

\end{thebibliography}
\end{document}